\documentclass[11pt]{amsart}

\usepackage{amsmath, amssymb,latexsym, amsthm, amsfonts, amscd}
\usepackage{verbatim}
\usepackage{epsfig}
\usepackage{graphics}

\title[Tangle Solutions for DNA-Rearranging Proteins]{Tangle Solutions for a Family of DNA-Rearranging Proteins
\mbox{\hspace{0.5in}\tiny{To Appear in} \hspace{0.5in}}
\mbox{\tiny\emph{Mathematical Proceedings of the Cambridge
Philosophical Society}}}

\author[D.~Buck]{Dorothy~Buck}
\address{Division of Applied Mathematics, Brown University}
\curraddr{Department of Mathematics, Imperial College London}
\email{d.buck@imperial.ac.uk}

\author[C.~Verjovsky~Marcotte]{Cynthia~Verjovsky~Marcotte}
\address{Department of Mathematics, St. Edward's University}
\email{cynthm@admin.stedwards.edu}

\newtheorem{thm}{Theorem}[section]
\newtheorem{THM*}{Theorem}

\newtheorem{lem}[thm]{Lemma}

\theoremstyle{definition}
\newtheorem{dfn}{definition}[section]
\theoremstyle{remark}

\numberwithin{equation}{section}

\newcommand{\nibf}{\noindent \textbf}

\newcommand{\ob}{O_c}
\newcommand{\of}{O_f^k}

\newcommand{\Qonly}{\mathbb{Q!}}
\newcommand{\wt}{\widetilde}
\newcommand{\bs}{\backslash}
\newcommand{\beq}{\begin{equation}}
\newcommand{\eeq}{\end{equation}}
\newcommand{\bea}{\begin{eqnarray}}
\newcommand{\eea}{\end{eqnarray}}

\newcommand{\dbc}{{\rm{dbc}}}
\newcommand{\core}{{\rm{core}}}
\newcommand{\interior}{{\rm{int}}}

\begin{document}

\maketitle

\section*{9 April 2004}

\medskip

\bigskip

\section*{Abstract}

We study two systems of tangle equations that arise when modeling
the action of the Integrase family of proteins on DNA.  These two
systems---direct and inverted repeats---correspond to two
different possibilities for the initial DNA sequence. We present
one new class of solutions to the tangle equations. In the case
of inverted repeats we prove that any solution not in these, or 2
previously known classes, would have to belong to one specific class.
In the case of direct repeats we prove that the three classes are
are the only solutions possible.

\section{Introduction}

We seek to mathematically model the action of a protein which
changes the topology of DNA. By  representing different regions
of a DNA molecule as tangles,  we can describe the protein's action
as a change in one of the constituent tangles.  We wish to find all
tangle combinations that explain the protein's action.

Ernst and Sumners \cite{ES1}, based on the biological work of
Wasserman and Cozzarelli \cite{WDC}, \cite{WC}, developed the
tangle model to describe and make predictions about how the
protein Tn3 interacts with DNA. The tangle model has since been
used to describe other protein-DNA interactions (see for example,
\cite{Cri}, \cite{Dar}, \cite{me}, \cite{SECS} and \cite{Vaz}).
Here, we focus on a particular family of proteins, the Integrase
family of recombinases (reviewed by Grainge and Jayaram
\cite{GJ}). Members of the Integrase family are involved in a wide
variety of biological activities, including integrating a virus'
DNA into a host cell's DNA (hence the name).

Integrase proteins share a common mechanism of rearranging DNA
sequences which can result in a change of DNA topology.
Specifically, by harnessing varying numbers of DNA axis
self-crossings (so-called \textit{supercoils}), these proteins
transform unknotted circular DNA into a variety of torus knots or
links \cite{GJ}.

For proteins such as Tn3, which require a fixed number of DNA
supercoils, before rearranging the DNA (often multiple times) and
then releasing it, it has been possible to completely solve the
model tangle equations (see \cite{ES1} and \cite{ES2}). Integrase
proteins, however, can act on DNA with any number of supercoils,
and act only once before releasing the DNA. The resulting DNA can
thus take a variety of knotted or linked forms. To accurately
model this varying supercoiling, and its effects on the resulting
DNA, we must use a larger number of tangles than needed for Tn3.
The increased complexity of the tangle model for the Integrase
family has thus far prevented a full solution to the tangle model.
Previous work, notably \cite{Cri} and \cite{me}, have found
solutions by making several assumptions which are thought to be
biologically reasonable and which lead to considerable
mathematical simplification.

The current work considers the tangle model, with no simplifying
assumptions, for a generic member of the Integrase family.   We
use the protein Flp as an illustration, but the results are
identical for other members of the Integrase family including
$\lambda$ Int and Cre.

We give three classes of solutions for each of the systems. Using
Dehn Surgery techniques we eliminate all possibilities (for
direct repeats) or all but one other possiblity (for inverted repeats).

As discussed in more detail below, we represent a circular DNA
molecule, before the protein action, by the numerator closure of
the sum of three biologically determined tangles: $\of, \ \ob$
and $P$.  We represent a DNA molecule after the protein
action as the numerator closure of the sum of $\of, \ \ob$ and
$R$.   We have two systems of tangle equations, corresponding to
two different initial types of DNA sequences:
\begin{tabbing}
\indent\= Before:\ \= \kill
\> \textrm{Before:\ } \> $N(\of + \ob + P)= b(1,1)$\\
\> \textrm{After:} \> $N(\of + \ob  + R) = b(2k,1)$ \ \textrm{or} \
$b(2k+1,1)$ { {\rule{0cm}{0.5cm}}}
\end{tabbing}
where $N(T)$ is the numerator closure of the tangle $T$, $b(p,q)$
is the $q/p$ four-plat, and $k \in \{0,1,2,3,4\}$. When the
products are of the form $b(2k,1)$ the system is called the
{\textit {direct}} case, and $b(2k+1,1)$ it is called the {\textit
{inverted}} case. We use $\{O_1, O_2\}=\{\of, \ob\}$;
\textit{i.e.}, the use of $O_1$ and $O_2$ means that we are
considering two cases: $O_1=\of$ and $O_2=\ob$, and viceversa. Our
main result is the following:

\begin{THM*}
There are three classes of solutions to the above equations, which are:
\begin{enumerate}
\item[1.] $P$ is the infinity tangle, $\ob$ and $\of$ are integral tangles
\item[2.] $P$ and $\ob$ are integral, $\of = (\infty)$ for at most
1 value of $k$, integral for at most 2 values of $k$ and otherwise
is a strictly rational tangle which is either vertical or the sum
of a vertical and a horizontal tangle.
\item[3.] $P$ is strictly rational,  $\ob$ is integral, and $\of$
is integral for at most 1 value of $k$ and strictly rational
otherwise.
\end{enumerate}

\vspace{-0.05in}

In the direct case, these are the only solutions.

In the inverted case, these are the only solutions with the
following possible exception: the double branched cover of
$\ob\cup\of$ is a hyperbolic (or trefoil for $k=2$) knot
complement, $P$ is rational, $O_1$ is integral and $O_2$ is
prime.  If $P$ is integral, then for at most one value of $k$,
$\of$ could be $(\infty)$.  If in addition, $\ob$ is also
integral then $\of$ could be integral for at most 2 values of
$k$, and  prime otherwise.
\end{THM*}

\vspace{-0.05in}

Theorem~1 is a consequence of the Theorems in Sections~6 and~7.  See
Table~1 for details.

From the biological point of view, the mechanism is expected to be
identical in both cases ({\textit{i.e.}}, $P$ and $R$\, in the
inverted case are the same as in the direct case) hence the
putative fourth possibility in the inverted case could be ruled
out as being biologically unfeasible. Whether there are
mathematical solutions to this class remains an open question.

The paper is organized as follows. In Section 2 we recall some
basic facts about tangles, four-plats and their corresponding
double branched covers. In Section 3 we provide the biological
motivation and background for our work: the action of the protein
Flp on DNA with either inverted or direct repeats. We also recall
Ernst and Sumners' mathematical model in terms of tangles and
four-plats \cite{ES1}. We present one novel class of solutions to
these tangle equations, and contrast them to two classes of known
solutions. In Section 4 we show that for direct repeats the
tangles $\of$ are rational, which extends the results of
\cite{ES1}. We also give an alternative proof to \cite{ES1} of the
fact that $R$ is rational. These proofs rely on Dehn Surgery
arguments for the complements of strongly invertible knots and for
Seifert Fibered Spaces. In Section 5, we prove a number of
technical results, principally concerning annuli in the double
branched covers of tangles. In Section 6, we use the results of
Section 5 to eliminate all but one of the remaining possibilities
for the tangles $\of$, $\ob$ and $P$. In Section 7, we place
further restrictions on the case that could not be entirely ruled
out for inverted repeats. We conclude with some remarks on the
biological relevance of our model and solutions, as well as
comments on further questions of interest. Table~1 can be used as
a guide to our results.

\section{Tangles, Four-Plats and their double branched Covers}

\begin{figure}[htb]
\begin{center}
\psfig{file=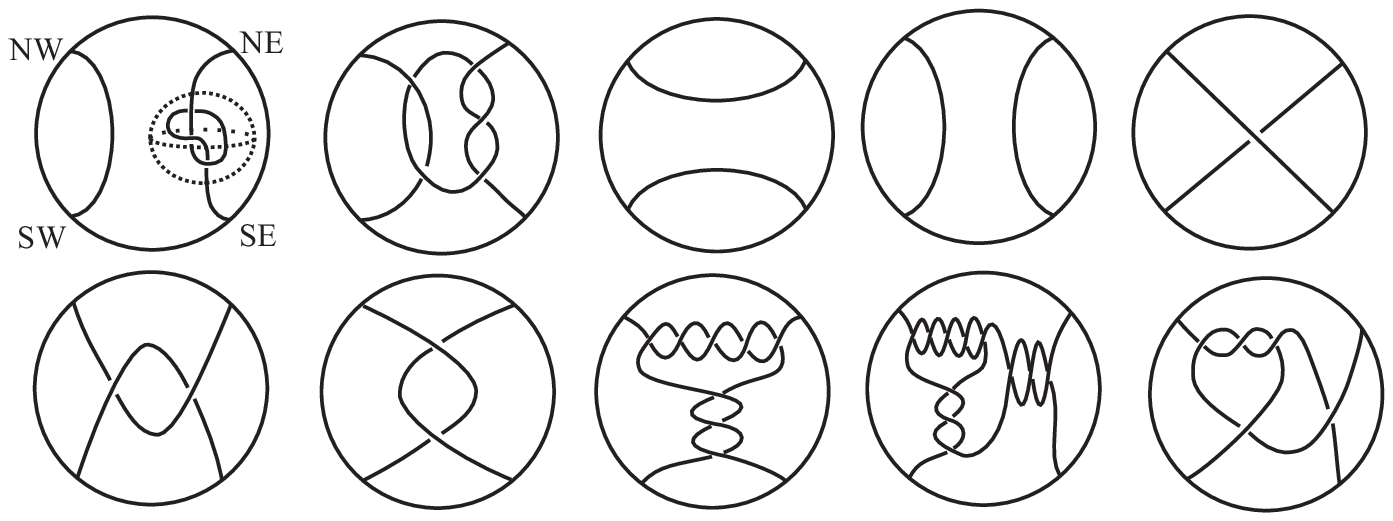,width=12.6cm}

\caption{Tangles. \mbox{{\em{Top Row:}} \/ Locally knotted, prime,
($0$), ($\infty$), ($1$)} \mbox{{\em{Bottom Row:}} \/ ($-2$),
($1/2$), ($5/14$), ($-37/14$), ($-7/4$)} \mbox{Note: The preferred
biological sign convention (shown above) is} opposite in sign to
that of Conway.} \label{f:tangleflavours}
\end{center}
\end{figure}

We begin by recalling a few elementary facts about tangles, and
the conventions that we will follow, established by Ernst and
Sumners in \cite{ES1}. A \textit{tangle}  $T$ is a pair
$(B^3,t)$, where $B^3$ is a 3-ball with a given boundary
parametrization with four distinguished boundary points labeled
NW, NE, SW, SE, and $t$ is a pair of properly embedded unoriented
arcs with endpoints NW, NE, SW and SE. See Figure
\ref{f:tangleflavours} for examples.

We say two tangles $A$ and $B$ are {\textit{equivalent}} if there
exists a series of moves that takes the strands of $A$ to
the strands of $B$, leaving the endpoints
fixed ({\textit{i.e.}, if there exists an isotopy from $A$ to
$B$ rel $\partial$).

Tangles can be divided into three mutually exclusive families:
locally knotted, prime and rational, illustrated in
Figure~\ref{f:tangleflavours}. A tangle is
\textit{locally knotted} if there exists a sphere in $B^3$
meeting $t$ transversely in 2 points such that the 2-ball bounded
by the sphere meets $t$ in a knotted spanning arc.  A tangle is
\textit{prime} if its double branched cover branched over $t$ is
irreducible and not a solid torus, as shown by Lickorish
\cite{Lic}.    All other tangles are \textit{rational}, so called
because their equivalence classes are in one-to-one
correspondence with the set of rational numbers (and infinity)
via a continued fraction expansion, as described by Conway
\cite{Con}.  A tangle whose corresponding rational number is
$p/q$ will be denoted by $(p/q)$. Rational tangles are formed by
an alternating series of horizontal and vertical half-twists of
two (initially untwisted) arcs. We say a rational tangle is
\textit{integral}, and write it as $(n)$, if it consists of a
series of $n$ horizontal half-twists, where $n \in \mathbb{Z}$.
We denote this class as $\mathbb{Z}$.  Similarly, a rational
tangle is \textit{vertical} if it consists of a series of $n$
 vertical half-twists ($|n|>1$), and denote it by $(1/n)$. We say
a tangle is  \textit{strictly rational}, and denote this class as
$\Qonly$ if it is neither integral nor the infinity tangle, and
so $\mathbb{Q} = \Qonly \cup \mathbb{Z} \cup \{ (\infty)\}$.

As Bleiler showed in \cite{Ble}, the minimal prime tangle has a
projection with five crossings (shown in
Figure~\ref{f:tangleflavours}).  Hence for certain tangles whose
strands represent very short DNA segments (\textit{e.g.}, less
than 50 base pairs), the most biologically relevant tangles are
those that are rational (See Maxwell and Bates \cite{MB} for more
details on the physical properties, including flexibility, of
DNA.)

There are several operations one can perform on tangles.  We
concentrate on three (see Figure~\ref{f:tangleops}).  Two of
these operations form a knot or link from a given tangle $A$: the
numerator closure, $N(A)$ and the denominator closure, $D(A)$.
The third operator, tangle summing, takes a pair of tangles $A$
and $B$  and (under certain restrictions) yields a third tangle,
$A +B$.  Note that the (0) tangle is the identity under this
operation: $A+(0)=A$. (See Conway \cite{Con}  for more details.)

\begin{figure}[htb]
\begin{center}
\psfig{file=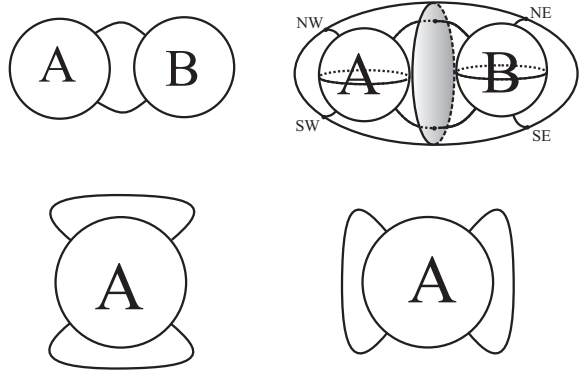,width=3in}
\caption{Three tangle operations. {\em Top:} tangle sum {\em Bottom:}
numerator closure and denominator closure}
\label{f:tangleops}
\end{center}
\end{figure}

The numerator or denominator closure of a rational tangle yields
a \textit{four-plat}, a knot or 2-component link that admits a
projection consisting of a braid on 4 strings, with one strand
free of crossings \cite{BS}   (see Figure \ref{f:4plat}).
Schubert showed that all four-plats are
prime knots and prime $2$-component links
(except for the unknot and the unlink of two unknotted components) \cite{Sch}.
A four-plat can be specified by a pair of integers $p$ and $q$,
and is written as $b(p,q)$. For example, we write the
unknot as $b(1,1)$, and the
trefoil as $b(3,1)$.

\begin{figure}[htb]
\begin{center}
\psfig{file=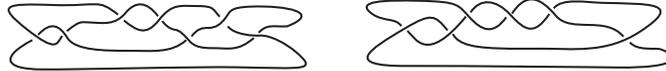,width=3.5in} \caption{The four-plats
$b(19,8)$ and $b(-9,4)$.} \label{f:4plat}
\end{center}
\end{figure}

If $T$ is a tangle, then we will write  $\wt{T}$ to  mean the
double cover of $B^3$, branched over the tangles arcs of $T$. We
will write $\dbc(K)$ to denote the three-manifold that is the
double  cover of $S^3$ branched over the set $K$. We now turn our
attention to the three-manifolds that arise as double branched
covers of tangles or four-plats. In this work all three-manifolds
are compact, connected and orientable.

If $P$ is a rational tangle, then $\wt{P}$ is a solid  torus,
which we will denote by $V_P$.   For notational simplicity we
will use $V^k_f$ and $V_b$ for the double brancheded covers of
$\of$ and $\ob$ respectively when they are rational tangles.
Schubert showed that $\dbc(b(p,q))$ is the  lens space $L(p,q)$.
Two four-plats $b(p,q)$ and $b(p',q')$ are equivalent iff their
corresponding double branched covers, the lens spaces $L(p,q)$ and
$L(p',q')$, are homeomorphic \cite{Sch}.  (See Rolfsen
\cite{Rolf} for more details on lens
spaces.)

The interplay of tangles and four-plats with their corresponding
double branched covers is the key to most of our results. When the
sum and subsequent numerator closure of two tangles $C$ and $D$
yields a four-plat $b(p,q)$, this induces a gluing of the
boundaries of their respective double branched covers $\wt{C}$ and
$\wt{D}$ that results in the lens space $L(p,q)$:
$$N(C+D) = b(p,q) \ \Leftrightarrow\  \wt{C} \cup_h  \wt{D} = L(p,q)$$
where $h:  \partial \wt{C} \rightarrow \partial \wt{D}$.
In particular, when $C$ and $D$ are both
rational, their  corresponding double branched covers are solid tori $V_C$
and $V_D$.  $V_C$ and $V_D$ form a Heegaard splitting of $L(p,q)$,
where $h$ is the map that takes $\mu_{\partial \wt{C}}$ to
$p\lambda_{\partial \wt{D}} + q\mu_{\partial \wt{D}}$.

\section{Biological Motivation, Model and Solutions}

We illustrate our model with the protein \textit{Flp} (pronounced
`flip'), a member the Integrase family of recombinases. Roughly
speaking, Flp recognizes two copies of a specific DNA sequence,
binds at these sites, cuts the DNA, moves the strands, reseals the
break, and releases the DNA. When acting on circular DNA, Flp can
change the underlying knot type of the DNA, for example turning
the unknot into the trefoil knot. The resulting knot type can
vary, depending on the amount of supercoiling of the DNA at the
time Flp binds. We call a DNA molecule that has not been acted on
by Flp a \textit{substrate}, and a molecule that has been acted on
a \textit{product}. In these terms then, the substrate is always
an unknot and the products are torus knots or links.


We model each of the substrates and products as the numerator
closure of the sum of three tangles, $N(\of + \ob + (P$ {or}
$R))$. Each tangle arc represents a segment of double-stranded
DNA. In the tangle model developed by Ernst and Sumners
\cite{ES1}, the cutting and joining of the DNA is assumed to be
completely localized: two of the tangles are unchanged by the
action of the protein. In the substrate, the first tangle, $P$
(Parental), represents two short identical sites that Flp
recognizes and to which it chemically binds and then cuts,
rearranges and re-seals. This action can be thought of as removing
$P$ and replacing it with a new tangle, $R$ (Recombinant) in the
product. The second tangle, $\ob$, represents the part of the DNA
that is physically constrained, but unchanged, by the protein ($O$
stands for Outside and $c$ for constrained). The last tangle,
$\of$, represents the part of the DNA that is free (hence the
subscript $f$) from protein binding constraints. $\of$ can vary
depending on the amount of DNA writhe present at the time Flp
acts. The superscript $k$ indexes these different possibilities.

Thus, the action of the protein is modeled by replacing $P$ in the
substrate by $R$ in the product (see Figure~\ref{f:recevent}). In
terms of tangles, this amounts to saying that our substrate and
products can be modeled as:

\begin{tabbing}
this is as long\= \kill
\>$N(\of + \ob + P) =$ \textrm{substrate (before recombination)}  \\
\>$N(\of + \ob + R) =$ \textrm{product (after recombination)}{ { \rule{0cm}{0.5cm}}}
\end{tabbing}
where $k \in \{0,1,2,3,4\}$. $\of$ varies as $k$ varies, so we
obtain different products, as described below. We use $O^k$ to
mean the part unchanged by Flp, that is, $\of\cup\ob$.
 We use $\{O_1, O_2\}=\{\of, \ob\}$;
\textit{i.e.}, the use of $O_1$ and $O_2$ means that we are
considering two cases: $O_1=\of$ and $O_2=\ob$, and viceversa.

The initial motivation for describing the action of Flp in terms
of tangles stemmed from a biological question:  When Flp binds to
DNA at $P$, how are the two sites oriented with respect to one
another? This was found to be antiparallel, by Grainge, Buck and
Jayaram \cite{me}, through a combination of mathematical proofs
and biological experiments.  In terms of tangles, experimental
work indicated that $\ob = (1)$.  Mathematical arguments, based
on a few additional biologically reasonable assumptions, then
concluded that $P = (\infty)$, and thus the sites must oriented
in antiparallel alignment.  Subsequent crystal structures by Chen
\textit{et al.} confirmed this alignment \cite{Rice}.


In the tangle model, the ultimate goal is to determine precisely
all tangles involved. Unlike in the paper described above or the
work of Crisona \textit{et al.}~\cite{Cri} or
Sumners \textit{et al.}~\cite{SECS}, where
biological considerations were taken into account, we will be
looking exclusively at what the mathematics alone can say.

\begin{figure}[htb]
\begin{center}
\psfig{file=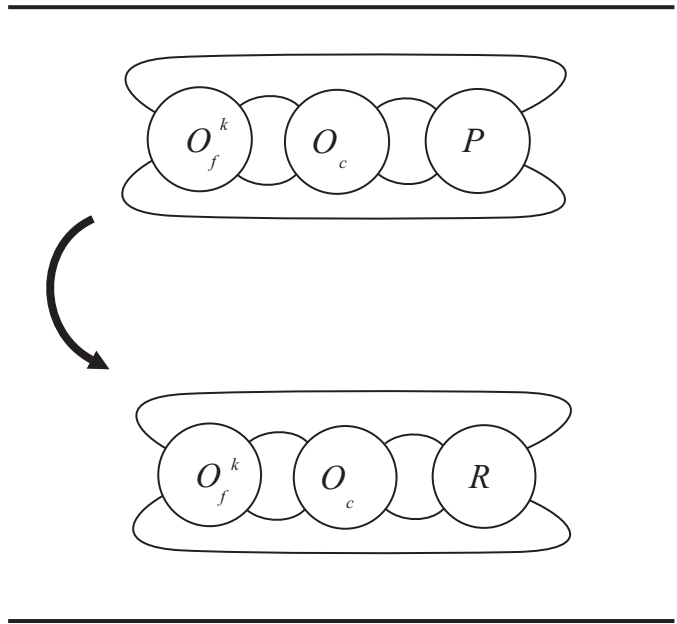,width=3.5in} \caption{Tangle Surgery}
\label{f:recevent}
\end{center}
\end{figure}

Replacing one tangle by another is known as tangle surgery.  If
each of $P$ and $R$ is rational, then tangle surgery corresponds
to replacing $V_P$ with $V_R$ in the double branched cover, and
thus corresponds to Dehn Surgery on $\wt{O^k}$.  Our  strategy is
to use restrictions on the type of Dehn Surgeries of $S^3=$
dbc$(b(1,1))$ that yield lens spaces.  This in turn restricts the
possible tangle solutions.

\subsection{Two Cases:  Inverted and Direct}

Flp identifies two short identical sequences, called
\textit{repeats}, on a molecule of DNA.  These sites have a
natural chemical orientation, and hence on circular DNA, the
strings can be in head to head (\textit{inverted repeats}) or
head to tail (\textit{direct repeats}) orientation
(Figure~\ref{f:dirinv}). When Flp acts on DNA it yields a variety
of knots or links that depend on $\of$. These products are torus
knots for inverted repeats, and torus links for direct repeats.
We consider each case separately.

\begin{figure}[htb]
\begin{center}
\psfig{file=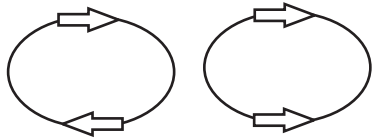} \caption{{\em Direct (left) and Inverted
Repeats.}} \label{f:dirinv}
\end{center}
\end{figure}

When Flp acts on a DNA molecule with inverted sites, experiments
have shown that the resulting DNA can be an unknot, or a knot
with up to 11 crossings \cite{GJ}.   Cozzarelli's lab has
obtained images (using electron microscopy) of the simplest
products, and has shown that they are the torus knots $b(1,1)$
(the unknot), $b(3,1)$ and $b(5,1)$ \cite{Cri}. This experimental
evidence indicates that Flp begins with an unknotted DNA
substrate with inverted repeats, $b(1,1)$ and converts it via
tangle surgery into a torus knot $b(2k+1,1)$, where $k \in
\{0,1,2,3,4\}$.   We thus model the action of Flp on DNA with
inverted repeats as:

\smallskip

\noindent\fbox{\begin{minipage}{12.4cm}
{\begin{center} \textsc{Inverted} \end{center}}
\begin{tabbing}
\label{invdown}
Before:\ \= \kill
\textrm{Before:} \>$N(\of + \ob + P) = b(1,1) =$ unknot, for \ $k \in \{0,1,2,3,4\}$\\
\textrm{After:} \>$N(O^0_f + \ob  + R) = b(1,1) =$ unknot { { \rule{0cm}{0.5cm}}}\\
\>$N(O^1_f + \ob + R) = b(3,1) =$ trefoil{ { \rule{0cm}{0.5cm}}}\\
\>$N(O^k_f + \ob + R) = b(2k+1,1) =$ {torus knot $T_{(2k+1,2)}$}
{ { \rule{0cm}{0.5cm}}}
\end{tabbing}
\end{minipage}}

\smallskip

When Flp acts on a DNA molecule with direct sites, experiments
have shown that the resulting DNA can be an unlink, or a
2-component link with up to 10 crossings \cite{GJ}.
Electrophoretic gels have determined that the simplest products
are $b(0,1)$, $b(2,1)$ and $b(4,1)$ \cite{GJ}. This experimental
evidence indicates that Flp begins with an unknotted DNA
substrate with direct repeats, $b(1,1)$ and converts it via
tangle surgery into a torus link $b(2k,1)$, where $k \in
\{0,1,2,3,4\}$.   We thus model the action of Flp on DNA with
direct repeats as:

\smallskip

\noindent\fbox{\begin{minipage}{12.4cm}
 {\begin{center} \textsc{Direct} \end{center}}
\begin{tabbing}
\label{invdown}
Before:\ \= \kill
Before: \>$N(\of + \ob + P) = b(1,1) =$ unknot, \textrm{for}\ $k \in \{0,1,2,3,4\}$ \\
After:\>$N(O^0_f + \ob + R) = b(0,1) =$ unlink { { \rule{0cm}{0.5cm}}}\\
 \>$N(O^1_f + \ob + R) = b(2,1) =$ Hopf link { { \rule{0cm}{0.5cm}}}\\
 \>$N(O^k_f + \ob + R) = b(2k,1) =$ {torus link $T_{(2k,1)}$}{ { \rule{0cm}{0.5cm}}}
\end{tabbing}
\end{minipage}}

\smallskip

\subsection{Solutions}

When proteins in the Resolvase family of recombinases, such as Tn3
recombine unknotted DNA, they first fix $\of := (0)$, then
rearrange the DNA, often multiple times, before releasing it. The
corresponding tangle equations:

\begin{eqnarray*}
\begin{array}{l}
N(O_f+O_b+P)=K_0  {\mbox{\hspace*{2.1cm} for the substrate}}\\
\hspace*{-1.5ex}\left. \begin{array}{l}
N(O_f+O_b+R)=K_1\\
N(O_f+O_b+R+R)=K_2 \\
\hspace*{2cm}\vdots\\
N(O_f+O_b+\underbrace{R+...+R}_{n})=K_n
\end{array}
\right\}{\mbox{for the products}}
\end{array}
\end{eqnarray*}

have been solved by Ernst and Sumners \cite{ES1}.

Integrase proteins, however, do not fix the DNA writhe and act
only once before releasing the DNA. The resulting DNA can thus
take a variety of knotted or linked forms. To accurately model
this varying writhe, and its effects on the resulting DNA
products, we must use a larger number of tangles than needed for
Tn3. We hence have not four tangles to find ($P$, $R$, $O_f$ and
$\ob$), but 3 fixed tangles ($P$, $R$ and $\ob$) together with the
family of tangles $\of$. This increased complexity of the tangle
model for the Integrase family has thus far prevented a full
solution to the tangle model.

Our ultimate aim, as mentioned above, is to describe all
solutions to the systems of  tangle equations that model the
action of Flp. For reasons that will become apparent, we
subdivide the possible tangles into the following categories:
locally knotted, prime, integral, the $(\infty)$ tangle, and
strictly rational. As we shall see, the first category is easy to
rule out, leaving the other four to consider.

We show there exists three classes of solutions.

Examples of one class of solutions was previously known
\cite{me}:  $P = (\infty)$ and the others are integral.  The
simplest example is when $P=(\infty)$, $R=(0)$ and $\of$ and
$\ob$ are such that $O^k=(\pm(2k+1))$ for inverted, or $O^k=(\pm
2k)$ for direct \cite{me}.

A second class is $P$ and $\ob$ are integral, and $\of = (\infty)$
for at most value of $k$, integral for at most 2 values of $k$ and
strictly rational otherwise. The simplest example in this class is
when $P$=(0), $R=(\infty)$ and $O^k = (\pm 1/(2k+1))$ for
inverted, or $O^0=(\infty)$ and $O^{k>0} = (\pm 1/2k)$ for direct
(see Figure~\ref{f:bioconj}), which is biologically equivalent to
the first example. (See the Conclusion for more details on
biological equivalence.) Using biologically reasonable arguments
that tightly constrain the mathematical possibilities (including
the following: since DNA is negatively supercoiled, then $\of =
(n)$ or $-n,0)$, and since $P$ and $R$ each represent short
segments of DNA, they have at most 1 crossing) Crisona \textit{et
al} rule out solutions of this type \cite{Cri}.

Another example, which is biologically non-equivalent, is when
$P=(2)$ and $R=(1)$. In this case for direct repeats $O^0_f=(-1)$,
$O^1_f=(-3)$ or $(-\frac{5}{3})$, $O^2_f=(-\frac{9}{5})$,
$O^3_f=(-\frac{13}{7})$ and $O^4_f = (-\frac{17}{9})$. For
inverted repeats $O^0_f=(\infty)$ or $(-\frac{3}{2})$,
$O^1_f=(-\frac{7}{4})$, $O^2_f=(-\frac{11}{6})$,
$O^3_f=(-\frac{15}{8})$ and $O^4_f = (-\frac{19}{10})$. See
Figure~\ref{f:AnothSol}.

In the second class $O_1$ in Table~1 must be $\of$.  Further, if $\of$
is strictly rational, it must be a tangle that can be written as
either vertical or the sum of a vertical and a horizontal tangle.
(Note that if the vertical and the horizontal twists have
opposite sign this is what is called a ``non-canonical" form.)

\begin{figure}[htb]
\begin{center}
\psfig{file=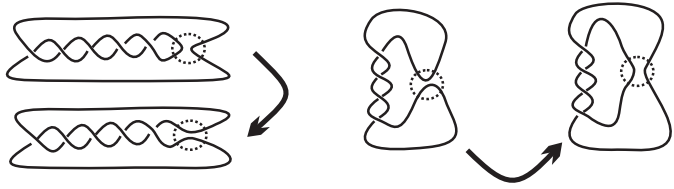,width=3.5in}
\caption{Biological Conjectures} \label{f:bioconj}
\end{center}
\end{figure}

In the following two theorems, we use the fact that
$N(A+n+\infty)=D(A)$.

\begin{thm}
\label{SecSol} If $P$ and $O_1$  are
integral and $O_2$ is strictly rational,
then in fact $O_1 = \ob$ and $O_2 = \of$.
In addition, $\of$ is either vertical, or the
sum of a vertical and a horizontal tangle. The horizontal
part must be equal to $-(P+\ob)$.
\end{thm}

\nibf{Proof.} If $P$ and $O_1$  are both integral then so is
$P+O_1$, \textit{i.e.,} $P+O_1=(m)$. As $O_2$ is strictly
rational, $O_2+(m)$ is also a strictly rational tangle. Since
$N(O_2+m)$ is the unknot, then $O_2+(m)$ must be vertical. This
means that $O_2= (1/n)+(-m)$ where $|n| >1$ (if $|n|=1$ then
$O_2$ is integral). The proof is complete, if $\of = O_2$.

Suppose instead $\ob=O_2=(1/n)+(-m), \ P=(s)$ and $\of
=O_1=(r_k)$.  Since $\of+\ob+P$ must be vertical, we have
$(1/n)+(-m)+(r_k)+(s)=(1/n)$.  Thus $(r_k)=(m-s)=\of$ is fixed,
which can happen for at most 1 value of $k$, say $k = j$.  For $k
\neq j$, if $\of = (\infty)$, then the substrate equation is
$N(\frac{1}{n} -m + s + \infty) = D(1/n) = b(1,1)$, and so
$|n|=1$.  This contradicts $|n| >1$.  Additionally, $\of$ cannot
be prime by Theorem~\ref{1int1prime}, or strictly rational by
Theorem~\ref{2strat}. Hence there are no possiblities for $\of$
when $k\neq j$.  Since the tangles equations for $k=j$ and $k
\neq j$  must be solved simultaneously, this case  cannot occur.
\hfill 

\begin{figure}[htb]
\begin{center}
\psfig{file=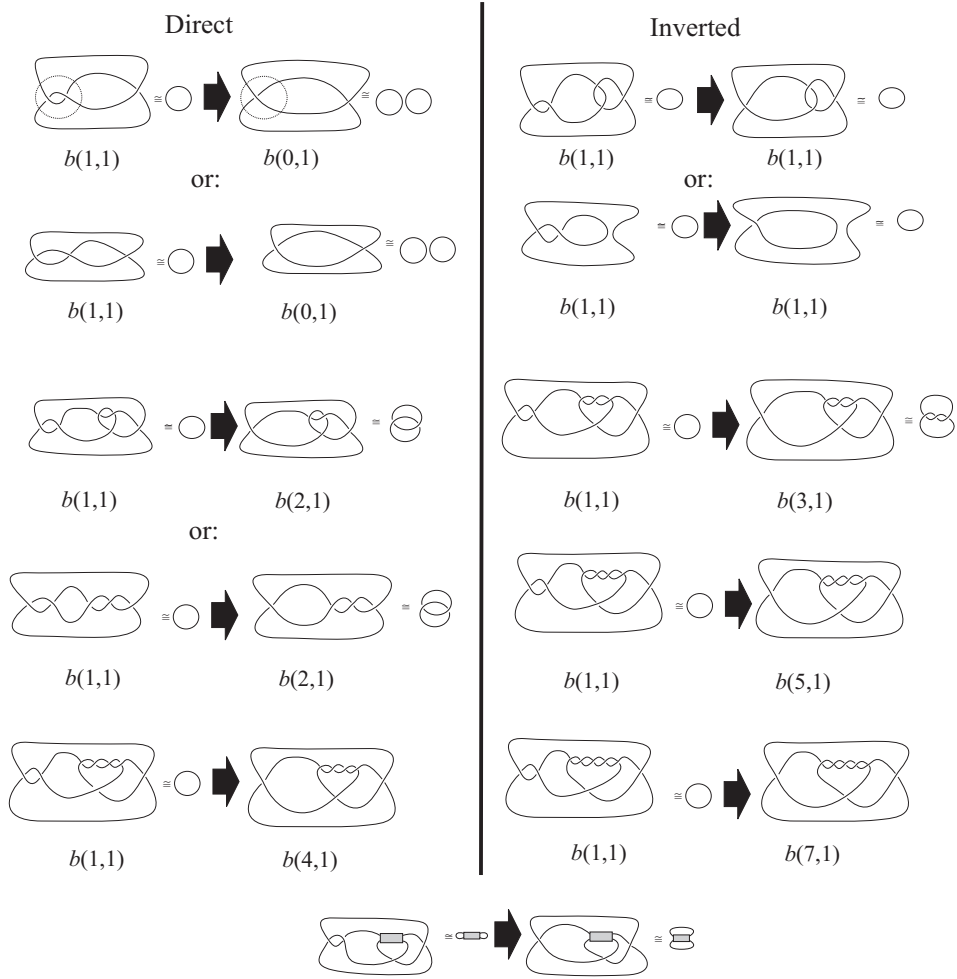,width=12.6cm} \caption{Solutions
when  $P=(2)$ and $R=(1)$. \emph{Direct repeats:} $O^0_f=(-1)$,
$O^1_f=(-\frac{5}{3})$ or $(-3)$, $O^2_f=(-\frac{9}{5})$ and
$O^3_f=(-\frac{13}{7})$ \emph{Inverted repeats:}
$O^0_f=(-\frac{3}{2})$ or $(\infty)$, $O^1_f=(-\frac{7}{4})$,
$O^2_f=(-\frac{11}{6})$, $O^3_f=(-\frac{15}{8})$ and
$O^4_f=(-\frac{19}{10})$. Note that the $O^k_f$'s are in
canonical form (all crossings of the same sign), and not in
vertical plus horizontal form.  The last figure shows how
$b(1,1)$ becomes $b(n,1)$ for $n$ horizontal twists.}
\label{f:AnothSol}
\end{center}
\end{figure}

\smallskip
In addition, in the second class, if a solution has
$\of=(\infty)$ for one value of $k$ and $\of$ integral for one
other value of $k$, then we can say several things about the
solution.

\begin{thm}
\label{ThirdSol}
Given a solution in which $\ob$ is integral, then
\begin{enumerate}
\item[($i$)] If $O^i_f=(\infty)$ for one value $i$ then $R$ must be vertical or $(\pm 1)$.
\item[($ii$)] If in addition $P$ is integral and $O^j_f$ is integral
for one or two other values of $k$, then $P\in\{0,\pm 2\}$.
\end{enumerate}
\end{thm}

\nibf{Proof.}

\begin{enumerate}
\item[($i$)] Suppose $\ob=(b)$ and $O^i_f=(\infty)$ for some value $i$. Then the
product is $N(\infty + b +R)=D(R)$, which must be a torus knot or
link.  As $R$ is rational by Theorem \ref{Rrat} for direct or
\cite{ES1} for inverted, this implies $R$ must be vertical,
$R=(1/r)$ with $|r| \geq 1$.
\item[($ii$)] Suppose $P=(p)$ and $O^j_f=(n_j)$ for one value
$j$, $j \neq i$. Then the substrate equation is
$N(n_j+b+p)=b(1,1)$ so $p=-(n_j+b)\pm 1$. The product is
$N(n_j+b+1/r)$, a torus knot or link. This means that $n_j+b=\pm
1$, hence $P=0$ or $\pm 2$. \hfill 
\end{enumerate}
\smallskip

The third and last class of solutions is where all tangles but
one are strictly rational and the remaining tangle is integral.
We note that if $P$ and $O_1$ are strictly rational, and $O_2$ is
integral,  then $\of$ must be the strictly rational tangle, $O_1$ as
shown in the following:

\begin{thm}
\label{lastone} If $P$ and $O_1$  are
strictly rational and $O_2$ is integral,
then in fact $O_1 = \of$ and $O_2 = \ob$.
\end{thm}

\nibf{Proof.} Assume instead that $\ob$ is strictly rational and
that $\of$ is integral for all $k$.  Then, by subsuming any
horizontal twists of $O^0_f$ into $\ob$, we can write $O^0_f =
(0)$ and $\of = (n_k)$ where $|n_k| >0$ for all $k \neq 0$.  So
$P = (x/y)$ and $\ob = (u/v)$, where $y,v
> 1$.

So $N(\ob+O^0_f+P)=N(u/v+ x/y) = b(1,1)$ and $N(\ob+\of+P)=N(u/v
+ {(u+n_{k}v)}/v) = b(1,1)$.  As in Ernst \cite{Ern}, we have
$|uy+xv| = 1$ and $|uy +xv +n_{k}vy| = 1$. Let $\sigma_1 =
\textrm{sign}(uy+xv)$,  and $\sigma_2 =
\textrm{sign}(uy+xv+n_{k}vy)$.

If $\sigma_1 = \sigma_2$, then $uy+xv = uy +xv+n_{k}vy$ so
$n_{k}vy = 0$.  This implies $n_k = 0$ for all $k$, contradicting
$|n_k|
> 0$ for $k \neq 0$.  So $\sigma_1 = \sigma_2$ for $k=0$ only.
For $k\neq 0$, $\of$ cannot be strictly rational nor
$(\infty)$  by Theorem~\ref{infty}(\textit{ii}), nor prime by
Theorem~\ref{infty}(\textit{iii}).  Since the tangle equations
for all $k$ must be solved simultaneously, this case
cannot happen.

If $\sigma_1 \neq \sigma_2$, then $-uy-xv = uy+xv+n_{k}vy$, so
$n_{k}vy = -2(uy+xv) = -2\sigma_1 |uy+xv| = -2\sigma_1$,
\textit{i.e.}, $n_{k}vy = \pm2$.  Recall, however, that $v$, $y >
1$, hence this is impossible. \hfill 
\smallskip

\begin{figure}[htb]
\begin{center}
\psfig{file=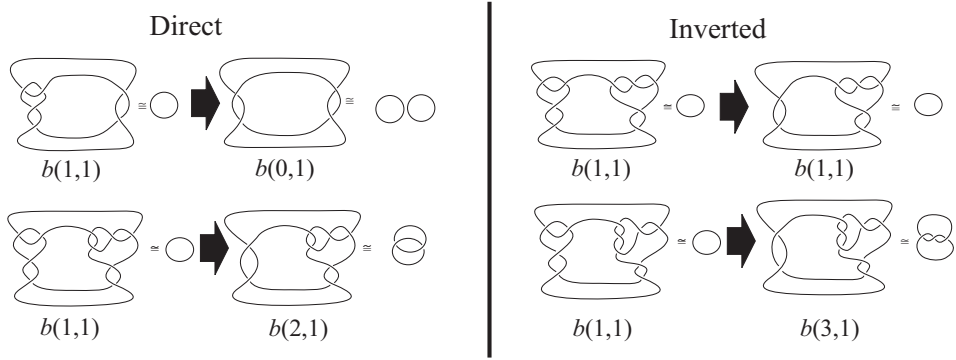, width=12.6cm}
\caption{Example of possible
solutions for $k= 0$ and $1$
when $P=(-\frac{2}{5})$.
In this case $R=(-\frac{1}{2})$.
In the direct case
$O^0=(\frac{1}{2})$ and $O^1=(\frac{5}{12})$
and in the inverted case
$O^0=(\frac{3}{7})$ and $O^1=(\frac{7}{17})$.}
\label{f:ThirdSol}
\end{center}
\end{figure}

\smallskip

An illustration of the last class of solutions is as follows:
Suppose $P$ is the rational tangle corresponding to some
non-integer rational, say, $(-\frac{2}{5})$ (see
Figure~\ref{f:ThirdSol}). Then tangles of the form
${U_n}=(\frac{1+2n}{3+5n})$ (for this value of $P$)
all satisfy $N(-\frac{2}{5}+{U_n})=$
unknot. Thus, all members of this class satisfy the first
condition (before recombination) in both the inverted and direct
cases. We must then find an $R$ such that for each $k$ there is a
$U_n$ such that $N(R+{U_n})= b(2k,1)$ (in the direct case) or $
=b(2k+1,1)$ (in the  inverted case). If $R=(-\frac{1}{2})$, then
the condition $N(R+{U_n})= b(\alpha, 1)$ is satisfied, where
$\alpha= 2(1+2n)-1(3+5n)=-(n+1)$. So for example, for $b(0,1)$ we
take $n=-1$, for $b(2,1)$ we take $n=-3$, and so on.

We now turn our attention to reducing the number of other
possibilities for the above tangles.
\section{Rationality of $\of$ for Direct Repeats}
We first note that all tangles in our model,
for both inverted and direct
repeats, must be locally \textit{un}knotted, since both the
substrate and one product are unknotted.  Therefore
all double branched covers of the tangles in our model are
irreducible, as are the double branched covers of their sums and
numerator closures.

Ernst and Sumners \cite{ES1} proved that $P$, $R$ and $O^0$
 (for both direct and inverted repeats) are rational tangles. In this section, we show
that $O^k = {O_f}^k + O_b$ rational for direct repeats, where $k
\in \{0,1,2,3,4\}$.
In order to prove $O^k$ is rational, we utilize the notion of
strongly invertible knots. A knot $K$ is \textit{strongly
invertible} if there exists an involution ({\textit{i.e.}}, an order 2
orientation-preserving homeomorphism), of $S^3$ which preserves
$K$ as a set and reverses the orientation of $K$. By a result of
Waldhausen, such an involution is equivalent to a $\pi$-rotation
whose axis is unknotted and which meets $K$ in exactly two points
\cite{Wal}.

We first establish some notation:  Let  $V_{(\infty)}$ be the
double branched cover of the $(\infty)$ tangle, and $r:
V_{(\infty)} \rightarrow V_{(\infty)}$ be the involution
corresponding to rotation by $\pi$ about a vertical axis that
punctures $V_{(\infty)}$ in two arcs.  The orbit space
$V_{(\infty)}/ r$ is a 3-ball, and the { induced} projection
$V_{(\infty)} \rightarrow V_{(\infty)}/ r$ is the branched
covering of ${(\infty)}$ branched over its arcs.

For any rational tangle $P$, and its solid torus double branched
cover $V_P$, there exist the following
two homeomorphisms: $m: {(\infty)} \rightarrow P$  which takes the
tangle arcs of ${(\infty)}$ to the tangle arcs of $P$, and
$\wt{m}: V_{(\infty)}\rightarrow V_P$ be such that the following
diagram commutes:
\begin{equation*}
\begin{CD}
r \mbox{\,{\resizebox{0.4cm}{!}{\includegraphics{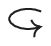}}}\,} V_{(\infty)} @> \wt{m} >> V_P \\
\indent @VVV        @V g VV\\
\indent {(\infty)}   @> m    >>  P
\end{CD}
\end{equation*}

\begin{lem}
\label{si}
Let $P$ be a rational tangle and $O^k$ a tangle such
that $N(O^k + P) = b(1,1)$ . Then $J={\rm{core}}(V_P)$ is a
strongly invertible knot in $S^3$.
\end{lem}

\nibf{Proof.} Recall that the double branched cover of $S^3$
branched over the unknot $b(1,1)$ is $S^3$.  Let $p$ be  the
covering map.  Then $p$ gives rise to an order two homeomorphism
$i:S^3 \rightarrow S^3$, where $i(x) = x$ if $x$ is in the branch
locus $p^{-1}(b(1,1))$, and $i(x) = y$ if $x \in S^3  \backslash
p^{-1}(b(1,1))$ and $\{x,y\} = p^{-1}(p(x))$.  Then $i$ is an
orientation preserving map and hence is an involution.  Note
$i|_{V_P} \cong h$, where $h$  is the homomorphism $\wt{m} \circ
r \circ \wt{m}^{-1}:V_P \rightarrow V_P$.  In fact,  $h$ is an
involution of $V_P$ such that $V_P / h$ induces the covering $g:
V_P \rightarrow P$ branched over  the two tangle arcs. Further,
we can choose $J' = \core(V_{(\infty)})$ so that it intersects
each branching arc exactly once and $r(J') = J'$.  Thus, by
choosing $J = \wt{m}(J) = \core(V_P)$, we ensure both that $J$
can be isotoped to intersect each branching arc in $V_P$ exactly
once and that $i(J)=J$. Since $r:V_{(\infty)} \rightarrow
V_{(\infty)}$ is a rotation by $\pi$, it reverses the orientation
on $J'$, and so $h = i|_{V_P}$ does so on $J$. Hence $J$ is a
strongly invertible knot in $S^3$.  \hfill 

\begin{thm}
$O^k = \of + \ob$ is rational in the case of direct repeats.
\end{thm}

\nibf{Proof.} Suppose $O^k$ is prime.  If $R$ is also prime,
then $\textrm{dbc}(N(O^k + R))$ must contain an
incompressible torus by Lickorish \cite{Lic},
but $\textrm{dbc}(N(O^k + R))=L(p,q)$, so we have a contradiction.

So if $\of$ is prime, then $R$ must be rational.  Now, by Lemma
\ref{si} $\wt{O^k_f}$ is the complement in $S^3$ of a
(nontrivial, as $O^k$ is prime) strongly invertible knot $J =$
\core$(V_P)$.  We then have $\wt{O^k} \cup V_P = L(1,1)$ for all
$k$, and $\wt{O^2} \cup V_R = L(4,1)$ and $\wt{O^3} \cup V_R =
L(6,1)$. This contradicts a result by Hirasawa and Shimokawa
\cite{HS}: there are no surgeries on nontrivial strongly
invertible knots that yield lens spaces of the form $L(2p,1)$ for
$p$ prime.  Hence, $\wt{O^k}$ must be the complement of the
trivial knot.  In other words $\wt{O^k}$ is a solid torus, and
thus $O^k$ is rational. \hfill 


\section{Background for Further Results}

Here we lay the necessary groundwork for Section 6, in which we
eliminate all but one class of possible solutions.

If $D$ is the gluing disc in the tangle sum $T = O^k + P$,
then it is straightforward to see that the $\wt{D} =
\dbc(D)$, branched over 2 points, is an annulus in $\partial\wt{P}
= \partial\wt{O^k}$.  We now consider $\partial\wt{D}$ in more
detail.

\begin{dfn}
A curve on the boundary of a solid torus will be said to
be \textbf{1-longitudinal, n-longitudinal}, or \textbf{meridional} if it
wraps once, $n$ times or $0$ times, respectively, longitudinally
around the solid torus boundary.
\end{dfn}

\begin{lem}
\label{forcing}
Given three rational tangles $R$, $S$ and $T$ (and gluing disc
$D$ between $S$ and $T$), such that $N(R+S+T) = b(1,1)$ and hence
$V_R \cup V_S \cup_{\partial{\wt{D}}} V_T = S^3$, then:
\begin{enumerate}
\item[($i$)] $\partial{\wt{D}}$ is meridional on $\partial{V_T}$ iff $T =
(\infty )$.
\item[($ii$)] $\partial{\wt{D}}$ is 1-longitudinal on $\partial{V_T}$ iff
$T$ is integral.
\item[($iii$)] $V_S \cup_{\wt{D}} V_T$ is a solid torus $W$ iff at
least one of $S$ and $T$ is integral.
\item[($iv$)] Given 3 rational tangles $R$, $S$, and $T$, such
that
 $N(R+S+T) = b(1,1)$, then one must be integral.
\item[($v$)] Given 3 rational tangles $R$, $S$, and $T$, such
that $N(R+S+T) = b(1,1)$, then if one is $(\infty)$
the other two must be integral
\end{enumerate}
\end{lem}

\nibf{Proof.}
\begin{enumerate}
\item[($i$ \& $ii$)]
It is simple to see that the lift of $D$, the gluing annulus for
the ($\infty$) tangle, is a meridional annulus on $\partial
V_{\infty}$. Every rational tangle can be obtained from the
($\infty$) tangle by an alternating series of vertical and
horizontal twists. These correspond in $\partial V_T$ to
longitudinal and meridional Dehn twists: a summing disc $D$ of
$(\infty )$ wraps equatorially around the 3-ball once for each of
the first vertical half-twists of $T$. In $V_T$, each vertical
half-twist of the strands of $T$  corresponds to one longitudinal
Dehn twist. The next horizontal  half-twists wrap $\wt{D}$
meridionally around $V_T$, which does not change the number of
longitudinal wraps.

Thus tangle $T$ can contain at most one vertical half-twist (plus
some number $n$ of horizontal half-twists, {\em{i.e.}}, $n+1$
horizontal half-twists) for $\partial{\wt{D}}$ be 1-longitudinal
on $\partial{V_T}$.  Therefore $T$ is integral.  The converse is
straightforward. (See Ernst \cite{Ern} for a complete description
of the construction of $\dbc(T)$.)

\item[($iii$)]  If $S = (\infty)$ then $\partial{\wt{D}}$
is meridional on $V_S$ and so bounds a merdional disc $\times I$.
Gluing this to a non-longitudinal curve would give a
punctured lens space $\not\cong S^3$, and so from ($ii$), $T$ must
be integral. If $S \neq (\infty)$,
gluing $V_S$ and $V_T$ together along annuli on their
boundaries yields a Seifert Fiber Space over a disc with 0, 1 or
2 exceptional fibers, where $\partial\wt{D}$ induces the
fibration. If $\partial{\wt{D}}$ wraps more than once
longitudinally on both $V_S$ and $V_T$, then the Seifert Fiber
Space has 2 exceptional fibers and hence cannot be a solid torus.
So $\partial{\wt{D}}$ is 1-longitudinal on at least one of
$\partial{V_S}$ or $\partial{V_T}$, and thus by ($i$) at least
one of $S$ and $T$ is integral.

Note that if $S$ (say) is integral, then by (\textit{ii}),
$\partial{\wt{D}}$ is 1-longitudinal on $\partial{V_S}$.  It is
straightforward to show that this implies that $V_S \cup V_T$ is
a solid torus.

\item[($iv$)]  If $V_R \cup V_S$ is a solid torus $W$, then by ($iii$)
$\partial{\wt{D}}$ is 1-longitudinal on $\partial{V_R}$ or
$\partial V_S$, {\em{i.e.}}, $R$ or $S$ is integral.  Otherwise,
suppose that $V_R \cup V_S = S^3\bs V_T$ is not a solid torus.
Then $S^3 \bs V_T = V_R \cup V_S$ is a non-trivial knot
complement, so $V_T$ is a torus knot neighborhood, and so $V_S
\cup V_T$ and $V_R \cup V_T$ are each solid tori.  Thus
$\partial{\wt{D}}$ is 1-longitudinal on $\partial{V_T}$, and so
by ($ii$) $T$ is integral.

\item[($v$)]   Assume one of the three tangles, say $R$, is
$(\infty)$. Then if $S$ and $T$ are strictly rational,  $N(S +
\infty + T) = D(S + T) = D(S) \# D(T) = b(1,1)$. This
means  $D(S)$ and $D(T)$ must each be the unknot.
Since for a rational tangle $A$  we have $D(A)$ is the unknot iff
$A$ is integral, $S$ and $T$ must both be integral.
\end{enumerate}

In Section 6, we also will need the following result on surfaces
in knot complements.

\begin{thm}[Simon \cite{Sim}]
\label{Sim}  If a knot $J$ in $S^3$ contains an essential annulus
in its complement, then $J$ is either a cable (possibly torus)
knot or composite knot.
\end{thm}

Let $X$ denote an essential annulus in $S^3\bs N(K)$.  In the
case of a composite knot, $X$  comes from a sphere $S$ that
intersects the knot twice. $S$ cuts $K$ into two arcs $\alpha$
and $\beta$, such that an arc joining the cut points on $S$
together with either $\alpha$ or $\beta$ will yield a nontrivial
knot.  In this case $\partial X$ is meridional on $N(K)$. If $K$
is a cable knot, then each component of $\partial{X}$ is
$n$-longitudinal along $N(K)$.

\section{Eliminating and Restricting Solutions}

We now have that $P$, $R$ and $O^0$ are rational for both direct
and inverted repeats, and that $O^k$ is rational for all $k$  for
direct repeats.  In the theorems that follow we examine different
possibilities for the tangles $\of$, $\ob$ and $P$, and either
eliminate or restrict each possibility in turn.  These theorems
hold for both direct and inverted repeats, unless explicitly
stated otherwise.

For each of the theorems below, we {use the fact  that $N(\of +
\ob + P) = b(1,1)$.  Recall we use $O_1$ and $O_2$ in the cases
where it is not necessary to distinguish between $\ob$ and
$\of$.  The first theorem follows from Lemma~\ref{forcing}.

\begin{thm}
\label{infty}{\rule{0cm}{0.1cm}}\newline
\vspace*{-0.5cm}
\begin{enumerate}
\item[($i$)] At most one of $\of$, $\ob$ and $P$ can be $(\infty)$.
\item[($ii$)] The three tangles $\of$, $\ob$ and $P$ cannot all be strictly
rational.  Also, the case $P$ and $O_1$ are strictly rational and
$O_2 = (\infty)$ cannot occur.
\item[($iii$)] The case $P = (\infty)$ or strictly rational, $O_1$
strictly rational  and $O_2$ prime cannot occur.
\item[($iv$)] If $P = (\infty)$ and both $\ob$ and $\of$ are rational,
then both $\ob$ and $\of$ must be integral.
\end{enumerate}
\end{thm}

\nibf{Proof.}
\mbox{\rule{0cm}{1cm}}
\begin{enumerate}
\item[($i$)] If there were two tangles that were
$(\infty)$, the circular substrate would include 2 closed
circles, or alternately would imply that $S^3$ contains a
punctured $S^2$ x $S^1$, a contradiction.

\item[($ii$)] This is a direct consequence of Lemma \ref{forcing}($iv$).

\item[($iii$)]  If $P$ and $O_1$ be rational, then $\wt{P}=V_P$ and
$\wt{O_1}=V_{O_1}$ and $V_{O_1}
\cup V_P$ must be a solid torus, otherwise $(\partial({V_{O_1}}
\cup V_P)$ would be incompressible.  Hence by Lemma
\ref{forcing}($iii$), either $P$ or $O_1$ is integral.

\item[($iv$)] Since $P = (\infty)$ this follows from Lemma~\ref{forcing}($v$). \hfill 
\end{enumerate}

\subsection{Restrictions on the knot core(\mathversion{bold}$V_P$\mathversion{normal})}

In the following three theorems, we show that, given certain
tangle types, $\core(V_P)$ must be either a composite, cable or
torus knot.  We then use Dehn surgery results to obtain a
contradiction.

\begin{thm}
\label{2primes} $\ob$ and $\of$ cannot both be prime.
\end{thm}

\nibf{Proof.} Suppose $\ob$ and $\of$ are prime---\textit{i.e.},
$\wt{\ob}$ and $\wt{\of}$ have incompressible boundaries. Since
$\ob$ and $\of$ are glued along a disc $D$, then $\wt{\ob}$ and
$\wt{\of}$ are glued along this disc's lift, the annulus
$\wt{D}$.  Note this annulus must be incompressible. So $S^3\bs
V_P$ is a knot complement containing an incompressible annulus.
Hence if $\partial{D}$ on $\partial{V_P}$ is meridional, $S^3\bs
V_P$ is composite, or if not  $S^3\bs V_P$ is cabled or torus (by
Theorem~\ref{Sim}).

If $V_P$ is a cable or torus knot, then since
$S^3\bs\{\wt{D} \cup V_P\} \cong  \interior(\wt{\of}) \cup
\interior(\wt{\ob})$, this
implies that either $\wt{\ob}$ or $\wt{\of}$ is a solid torus whose core
is the companion of the cable knot.  Thus either $\ob$ or $\of$
must be rational, contradicting our assumption.

Hence $S^3\bs V_P = \wt{\of} \ \cup  \ \wt{\ob}$ is a composite
knot complement.  This contradicts  a result of Bing and Martin
\cite{BingM}:  surgery on a composite knot never yields a lens
space.  \hfill 

\begin{thm}
\label{1int1prime}
 The case where $P$ is integral, $O_1$ is strictly rational and $O_2$ is
prime cannot occur.
\end{thm}

\nibf{Proof.}
Given $P$ is integral, then by Lemma \ref{forcing} (\textit{iii}),
$W = V_{O_1}\cup_{\wt{D}} V_P$
is a solid torus.  Let $D$ be the gluing disc in $P + O_1$,
then $\wt{D}$ is 1-longitudinal on
$\partial{V_P}$ and $n$-longitudinal ($n>1$) on
$\partial{V_{O_1}}$  by Lemma~\ref{forcing} (\textit{ii}).
Hence $\core(V_P)$ is isotopic to an $n$-longitudinal
curve on $\partial{V_{O_1}}$ and $\core(W)$ and
$\core(V_{O_1})$ are isotopic. Thus, $\core(V_P)$ is a $(m,n)$
cable of $\core(V_{O_1})$.
 Therefore $\wt{O_2} \cong S^3 \bs (V_{O_1} \cup
V_P) \cong S^3 \bs V_{O_1}$, \textit{i.e.}, $\wt{\of}$ is the
(nontrivial, as $O_2$ is prime) complement of the companion
knot $\core(V_{O_1})$ of the cable knot
$\core(V_P)$.

$S^3 \bs V_P = V_{O_1} \cup \wt{O_2}$ is the complement of a
nontrivial cable knot.  Removing $V_P$ and replacing it with
$V_R$ to obtain a new four-plat $b(p,q)$ is equivalent to surgery
on $\wt{O_2} \cup V_{O_1}$ which yields a lens space $L(p,q)$.
By a result of Bleiler and Litherland \cite{BL}, $p$ must be
greater than or equal to 23, which contradicts the fact that all
the products have $p\leq 7$. \hfill 

\begin{thm}
\label{2strat}
 The case where $P$ is integral, and both $\ob$ and
$\of$ are strictly rational cannot occur.
\end{thm}

\nibf{Proof.}
Since all three tangles are rational, we have $V^k_f
\cup_{\wt{D}} V_b \cup_{\wt{A}} V_P$, where $\partial{\wt{A}}
= \partial{\wt{D}}$.  Since $P$ is integral, $\partial{A}$ is
1-longitudinal on $\partial{V_P}$ (by Theorem~\ref{forcing}({\textit{ii})),
 $s$-longitudinal on $\partial{V_b}$
( with $s>1$) and $t$-longitudinal on $\partial{V^k_f}$ (with
$t>1$). Let $W = V_P \cup V_b$, a solid torus by
Lemma~\ref{forcing}(\textit{iii}).  Then $W \cup V^k_f$ is a
Heegaard splitting of $S^3$, and $\core(V_P)$ is isotopic to an
$s$-longitudinal curve  on $\partial{W}$.  Similarly, $W' = V_P
\cup V^k_f$ is a solid torus, with $\core(V_P)$ $t$-longitudinal
on $\partial{W'}$.  Thus $\core(V_P)$ is an $(s,t)$ torus knot,
such that surgery on its complement yields in particular $L(2,1)$
for direct or $L(3,1)$ for inverted repeats.  This contradicts a
corollary from Moser \cite{Mos}:  if surgery on a torus knot
complement yield a lens spaces $L(p,q)$, then $p \geq 5$. \hfill

\begin{thm}
\label{manycases}
 \mbox{\rule{0cm}{1cm}}
\begin{enumerate}
\item[1.]  For direct repeats, the following cases cannot occur:
\begin{enumerate}
\item[($i$)] $O_1$ is integral and $O_2$ is prime.
\item[($ii$)] $O_1$ and $O_2$ are both strictly rational.
\item[($iii$)]  $O_1$ is strictly rational and $O_2 = (\infty)$.
\end{enumerate}

\item [2.]  For inverted repeats, we have the following
restrictions:
\begin{enumerate}
\item[($i$)] if $\ob$ is integral then $O^0_f$ cannot be prime
\item[($ii$)]  if $\ob$ is strictly rational then $O^0_f$ must be integral
\item[($iii$)] if $\ob$ is prime, then $O^0_f = (\infty)$ \hfill 
\end{enumerate}
\end{enumerate}
\end{thm}

\nibf{Proof.} We have shown in previous sections, that $\ob +
\of$ is rational for direct repeats, and $O^0 = O_f^0 + \ob$ is
rational for inverted.  The result then follows directly from a
result of Quach \cite{Qua}: the sum of 2 tangles is prime if
either both subtangles are prime, one subtangle is prime and the
other rational ($\neq (\infty)$), or both are non-integral
rational.  \hfill 

\subsection{Restrictions from the range of \mathversion{bold}$k$\mathversion{normal}}

We will need this lemma for the following theorem.

\begin{lem}
\label{qonlyint} Suppose one of $P$, $\ob$, $\of$  is strictly
rational and the other two are integral.  Then the strictly
rational tangle can be written as
a rational tangle of length two, {\em
i.e.}, as one set of vertical twists followed by one set of
(possibly zero) horizontal twists.
\end{lem}

\nibf{Proof.} Since a strictly rational tangle plus an integral
tangle yields a rational tangle, we can write $N(A+B) = b(1,1)$,
such that, without loss of generality, $A = (a/b)$ and $B = (n)$.
Then by a result of Ernst and Sumners \cite{ES1}, $|a+bn| = 1$,
and so $\frac{a}{b} = \mp n + {1/{\pm b}}$, so $A$ is composed of
$\pm b$ vertical twists and $\mp n$ horizontal twists. \hfill

Note that the length two projection of the tangle may not the
canonical form, in which the projection of a rational tangle has
all the crossings of the same sign.In other words, a
non-canonical form occurs when the horizontal and the vertical
twists are of opposite sign.

The next theorem relies on the fact that for each $k \in
\{0,1,2,3,4\}$, $\of$ must be distinct.

\begin{thm}
\label{moreint}We have the following restrictions on $\of$:

\begin{enumerate}
\item[($i$)] If $P$ is strictly rational and $\ob$ is integral,
then $\of$ is integral for at most 1 value of $k$.

\item[($ii$)] If $P$ and $\ob$ are integral, then $\of$ is
integral for at most 2 values of $k$.
\end{enumerate}
\end{thm}

\nibf{Proof.}

\begin{enumerate}
\item[($i$)] If $P$ is strictly rational, then by Lemma
\ref{qonlyint}, $P = (p + 1/q)$,  $\ob = (b)$ and $\of = (f_k)$.
Since $N(f_k+b+ p\pm1/q)=b(1,1)$, then $f_k = -(b+p)$, a constant
independent of $k$.  This can occur for at most 1 value of $k$,
since if $i \neq j$, $O^i_f \neq O^j_f$.

\item[($ii$)] Similarly, if $\of = (f_k)$, $\ob = (b)$ and $P = (p)$,
then $N(\of + \ob +P) = N(f_k + b + p) = b(1,1)$.  So $f_k =
(-(b+p) \pm 1)$, \textit{i.e.}, there are at most two choices for
$f_k$.   \hfill 
\end{enumerate}

The following theorem is based on the fact that the denominator
or numerator closure of a single tangle cannot give rise to a
spectrum of products. We also use the following facts:
$N(A+\infty)=D(A)$, $D(A+n)=D(A)$ where $(n)$ is an integral
tangle, and $D(A+B) = D(A) \# D(B)$.  For $S$ rational,
$N(S)=b(1,1)$ iff $S$ is vertical, and $D(S)=b(1,1)$ iff $S$ is
integral.

\begin{thm}
\label{moreinfty}
The following cases cannot occur:
\begin{enumerate}
\item[($i$)] $P$ strictly rational, $O_1 = (\infty)$ and $O_2$ integral.
\item[($ii$)] $P$ integral, $O_1 = (\infty)$ and $O_2$ strictly rational.
\item[($iii$)] $P$ strictly rational, $O_1 = (\infty)$ and $O_2$ prime.
\item[($iv$)] $P$ integral, $\ob = (\infty)$ and $\of$ integral.
\item[($v$)] $P$ integral, $\ob =(\infty)$ and $\of$ prime.
\end{enumerate}
In addition, for direct repeats:
\begin{enumerate}
\item[($vi$)] The case $P$ integral, $\ob$ prime and $\of=(\infty)$  cannot occur.
\item[($vii$)] If $P$ and $\ob$ are integral and $\of =(\infty)$
for a value of $k$, then $\of$ may be integral for up to 2
values of $k$ and must be strictly rational for all
other values of $k$.
\end{enumerate}
For inverted repeats:
\begin{enumerate}
\item[($viii$)] If $P$ is integral, $\ob$ is prime, and $\of =(\infty)$
for a value of $k$, then $\of$ must be integral for all
other values of $k$.
\end{enumerate}
\end{thm}

\nibf{Proof.}
\begin{list}{($i$)}{\labelsep .2cm \itemindent 1.3cm \labelwidth 2cm \leftmargin 1.4cm }
\item[($i$, $ii$ \& $iii$)]
Note first that the substrate equation for either inverted or
direct is $N(\infty + O_2 + P) = D(O_2 +P) = D(O_2) \ \# \ D(P) =
b(1,1)$.  Thus both $D(O_2)$ and $D(P)$ must be the unknot.
Hence $O_2$ must be integral or prime. In addition, as $P$ is
rational, $P$ must be integral. Therefore Cases ($i$), ($ii$) and
($iii$) are impossible.
\end{list}
\begin{enumerate}
\item[($iv$)]
If $\ob = (\infty)$, and $\of = (n_k)$,  then the product is
$N(n_k + \infty + R) = D(R)$, which is a unique four-plat.  Hence
$\of$ is integral for at most one value of $k$, say $k=j$.   For
$k \neq j$, then as mentioned above, $\of$ must be prime.  The
product equations are: $N(\infty + \of + R) =D(\of) \# D(R) =
b(2k(+1),1)$.  As noted above, given the substrate equation
$D(\of)$ must be the unknot, and thus the product is $D(R)$,
which is a fixed four-plat. Hence there can be only one value of
$k$, say $k = i$, such that $O^i_f$ is prime. Therefore this case
cannot occur since this accounts for only 2 values of $k$, but $k
\in \{0,1,2,3,4\}$.

\item[($v$)]
If $\ob = (\infty)$, and $\of$ is prime, then, arguing as in
Case~($iv$), $\of$ can be prime for at most one value of $k$, say
$k=j$.   For $k\neq j$, then as shown in Cases~($i$), ($ii$) and
($iii$), $\of$ must be integral.  As in Case~($iv$), $\of$ can be
integral for at most 1 value of $k$ as well, and so this case, as
above, cannot occur.

\item[$(vi)$]
In the case of direct repeats,
suppose $P$ is integral
and $\ob$ is prime, and for some
$j$, $O^j_f = (\infty)$.  Then for $k \neq j$, $\of$ cannot be
prime, by Theorem~\ref{2primes}.  In addition
$\of$ cannot be strictly rational, by Theorem~\ref{1int1prime}.
Finally $\of$ cannot be integral
by Theorem~\ref{manycases}($i$), and so this case cannot
occur.

\item[($vii$)]
In the case of direct repeats,
if $P$ and $\ob$ are integral and for some $j$, $O^j_f =
(\infty)$, then for $k\neq j$, $\of \neq (\infty)$.  If $\of$ is
integral, then there can be at most two values of $k$ such that
$\of$ is integral by  Theorem~\ref{moreint}($ii$).
For direct repeats, $\of$ cannot be prime by
Theorem~\ref{manycases}($i$). Thus $\of$ must be
strictly rational for all other values of $k$.  Note that if
$\of$ is strictly rational then by Theorem~\ref{SecSol} $\of$
must be either a vertical tangle or the sum of a vertical and a
horizontal tangle.

\item[($viii$)]
In the case of inverted repeats, suppose $P$ is integral and
$\ob$ is prime, and for some $j$, $O^j_f = (\infty)$.  Then for
$k \neq j$, $\of$ cannot be prime, by Theorem~\ref{2primes}.  In
addition $\of$ cannot be strictly rational, by
Theorem~\ref{1int1prime}.  Therefore, for $k \neq j$, $\of$ must
be integral. \hfill  
\end{enumerate}
\section{Restrictions on Remaining Inverted Case}

In Section~6 we eliminate all classes of solutions, other than
the three classes given in Theorem~1, for direct repeats.

In the inverted case, these three classes of solutions are the
only ones, with possible exception of the case where $P$ is
rational, $O_1$ is integral and $O_2$ is prime.  If $P$ is
integral, then for at most one value of $k$, $\of$ could be
$(\infty)$.  If in addition, $\ob$ is also integral then $\of$
could be integral for at most 2 values of $k$, and strictly
rational or prime otherwise. We now prove that $\wt{O^k}$ is the
complement of a hyperbolic knot or, for $k=2$, possibly the
trefoil knot.



\begin{thm}
\label{torsat} For inverted repeats, if $P$ is strictly
rational,  $O_1$ is integral and $O_2$ is prime,
then  $\wt{O^k}$ cannot be a satellite knot complement.
For $k\neq 2$ $\wt{O^k}$ cannot be a torus knot complement
For $k = 2$, then if $\wt{O^k}$ is a torus knot, it
must be the trefoil knot complement.
\end{thm}

\nibf{Proof.} \ \textit{Case 1:} $\core(V_P)$ is a satellite
knot.  For both inverted and direct repeats this contradicts a
result of Bleiler and Litherland \cite{BL}: if surgery on
$\core(V_P)$ yields a lens space $L(p,q)$, then $p \geq 23$.

\textit{Case 2:} $\core(V_P)$ is a (non-trivial) torus knot,
$T_{(m,n)}$, so $(m,n) = 1$, and as it is non-trivial $m, n>1$
and $mn\geq 6$.  If $T_{(m,n)}(r)$, the manifold resulting from
surgery along slope $r$, is a lens space,  results of Moser
\cite{Mos} tell us that the lens space must be $L(p,qn^2)$ where
$p$ and $q$ are related to $m$ and $n$ by $p = qmn \pm 1$.  Since
for $T_{(m,n)}$ nontrivial, $mn \geq 6$, it follows that $p \geq
5$. It is simple to check that the smallest lens spaces that can
result are $L(5,4)$, $L(7,2)$, $L(7,4)$, $L(9,4)$ and $L(9,7)$.
By Brody \cite{Bro}, only $L(5,4)$ is homeomorphic to a product,
namely $L(5,1)$ (the product when $k=2$). $L(5,1)$ can only be
obtained by surgery on the trefoil knot complement. For the other
products $L(2k+1,1)$ ($k \neq 2$), we obtain a contradiction.
\hfill 

\smallskip

A simple closed curve on a handlebody $W$ is \textit{primitive}
if it is a free generator of $\pi_1(W)$.  A curve on $F=W_1\cap
W_2$, a Heegaard splitting surface,
 is \textit{doubly primitive} if it is primitive on both handlebodies $W_1$
and $W_2$.  A \textit{Berge knot}, K, is a doubly primitive curve
on a genus two Heegaard surface of $S^3$. Berge~\cite{Ber} showed
that there is an integer surgery on $K$ that yields a lens space.
The Berge Conjecture states that this is the only type of
non-trivial knots that have surgery yielding lens spaces.

In the current setting, from Theorem \ref{torsat}, for $k\neq 2$,
$\wt{O^k}$ would have to be
the complement of either a hyperbolic Berge knot
or a knot that is a counterexample to the Berge conjecture.

\section{Conclusion and Directions for Further Research}

We have given three families of solutions for the tangle
equations that arise from the action of the protein Flp on DNA.

Recall in Section 3.2 we gave three classes of solutions.  The
second class  included the example where $\of$ is  strictly
rational, $\ob$ is integral, $P=(0)$ and $R=(\infty)$.   This
example, however, is \textit{biologically equivalent} to the
solution introduced in \cite{me}:  $P$ is the infinity tangle,
$\ob$ and $\of$ are integral tangles.  That is to say, the two
classes of solutions reflect identical models of protein
action---tangle surgery interchanges the only two tangles with
zero crossings  ($(\infty)$ and $(0)$). The second example given
of this class, however, is not biologically equivalent. The third
class of solutions includes many different (biologically
non-equivalent) possible actions.

We have eliminated all other classes of solutions in the case of
direct repeats, and all but one in the case of inverted. Several
open questions remain.  The first two questions are directed
towards biologists.  Firstly, although experiments have
determined the number of crossings of the product knots and
links, for higher crossing products it has not been confirmed
that they are in fact torus knots and links. Can we obtain
experimental confirmation of this? Secondly, although the third
class of solutions is mathematically possible, biological
considerations such as DNA's stiffness impeding a high number of
crossings, make them biologically unlikely.
So can the third class of solutions be ruled out experimentally?

The third question is mathematical in nature, and is as follows:
as noted in Section 7, one
possibility remains open for the inverted case:  $ P\in
\mathbb{Q}$,  $O_1\in \mathbb{Z}$, and $O_2$ prime.  It is
therefore natural to ask whether this possibility can be
eliminated. We plan to investigate whether there exist examples
of $\wt{O^k}$ that are the complement of a hyperbolic Berge knot
or (for $k=2$) the trefoil knot. In addition, it may be possible
to eliminate---or find solutions within---certain subclasses of
prime tangles, such as Montesinos knots.

Again, although we have described the biology and tangle
solutions in detail for the protein Flp, our results hold more
generally.  In particular, they can be applied to all proteins in
the Integrase family, such as $\lambda$ Int and Cre, whose
products are torus knots or links.

\smallskip

\nibf{Acknowledgements.} We wish to thank Cameron Gordon and John
Luecke for a number of illuminating discussions. We also would
like to thank Makkuni Jayaram for introducing us to Flp. Finally,
we wish to thank Mark Haskins for his careful reading of previous
versions, and the reviewer for his/her careful reading and
insightful comments. DB was supported by a grant from the
Burroughs Wellcome Fund. CVM was supported in part by the
Presidential Excellence Summer Scholarly Activity Grant from St.\
Edward's University.

\bibliographystyle{amsplain}

\newpage

\begin{table}
    \begin{tabular}{||c|c|c|l|l||} \hline \hline

        $P$ & $O_1$ & $O_2$ & Theorem where case & \\
             &       &      & is eliminated or restricted & Solutions\\
  \hline \hline
        $\infty$ & prime & prime & Theorem \ref{2primes} & \\
  \hline
                  &  $\Qonly$ &prime & Theorem \ref{infty}({\it{iii}})& \\
  \hline
                  & & $\mathbb{Z}$ & Theorem \ref{infty}({\it{iv}}) & \\
  \hline
                    &  & $\Qonly$ & Theorem \ref{infty}({\it{iv}}) & \\
  \hline
                   & $\mathbb{Z}$  & prime    & Theorem \ref{manycases}({\it{i}}) for direct & \\
                     &  &       &Theorem \ref{torsat} for inverted,  & \\
                    &           &          & {\hspace*{.5cm}}when $\wt{O^k}$ is a torus ($k\neq2$)  & \\
                     &           &          & {\hspace*{.5cm}}or satellite knot complement & \\
  \hline
                    &   & $\mathbb{Z}$  &   & SOLUTION 1\\
  \hline \hline

        $\mathbb{Z}$   & prime  & prime   & Theorem \ref{2primes} & \\
  \hline
                    & $\infty$ & prime & Theorem \ref{moreinfty}({\it{v}}) $\ob=O_1$ & \\
                    &           &       &Theorem \ref{moreinfty}({\it{vi}}) $\ob=O_2$ for direct &  \\
  \hline
                    &  & $\mathbb{Z}$ & Theorem \ref{moreinfty}({\it{iv}}) $\ob=O_1$ & SOLUTION 2  \\
  \hline
                    &  & $\Qonly$ & Theorem \ref{moreinfty}({\it{ii}}) & \\
  \hline
                    &$\Qonly$  & prime   & Theorem \ref{1int1prime} & \\
  \hline
                    & & $\Qonly$ &  Theorem \ref{2strat} & \\
  \hline
                    &$\mathbb{Z}$ & $\mathbb{Z}$ & Theorem \ref{moreint}({\it{ii}})
                                            for all but 2 values of $k$ & SOLUTION 2\\
  \hline
                    & & prime & Theorem \ref{manycases}({\it{i}}) for direct  &\\
                     &  &       &Theorem \ref{torsat} for inverted,  & \\
                    &           &          & {\hspace*{.5cm}}when $\wt{O^k}$ is a torus ($k\neq2$)  & \\
                     &           &          & {\hspace*{.5cm}}or satellite knot complement & \\
  \hline
                    & & $\Qonly$ & $\ob = O_1$ by Theorem \ref{SecSol} & SOLUTION 2\\
  \hline \hline

        $\Qonly$    & prime     &prime   & Theorem \ref{2primes} & \\
  \hline
                    & $\infty$  &prime & Theorem \ref{moreinfty}({\it{iii}}) & \\
  \hline
                    &           & $\mathbb{Z}$ & Theorem \ref{moreinfty}({\it{i}}) & \\
  \hline
                    &           & $\Qonly$ & Theorem \ref{infty}({\it{ii}}) & \\
  \hline
                    & $\Qonly$ & prime & Theorem \ref{infty}({\it{iii}}) & \\
  \hline
                    &          &$\mathbb{Z}$ & $O_2=\ob$ by Theorem \ref{lastone} &  SOLUTION  3\\
  \hline
                    &          &$\Qonly$   & Theorem \ref{infty}({\it{ii}}) & \\
  \hline
                    & $\mathbb{Z}$  & $\mathbb{Z}$ & Theorem \ref{moreint}({\it{i}}) for all but 1 value of $k$ & SOLUTION 3 \\
  \hline
                    &           & prime  & Theorem \ref{manycases}({\it{i}}) for direct & \\
                    &  &       &Theorem \ref{torsat} for inverted,  & \\
                    &           &          & {\hspace*{.5cm}}when $\wt{O^k}$ is a torus ($k\neq2$)  & \\
                     &           &          & {\hspace*{.5cm}}or satellite knot complement & \\

  \hline \hline
\end{tabular}

\medskip

\caption{All possible tangles classes for $P$ and $\{O_1, O_2\}=\{\ob,\of\}$.
$P$ is rational \cite{ES1}, and at most one of $P$,
$O_1$ and $O_2$ is $(\infty)$ by Theorem \ref{infty}.}

\end{table}

\clearpage

\appendix
\section{$R$ is rational for Direct Repeats}

We provide an alternative proof to Ernst and Sumners' result
\cite{ES1} that $R$ must be rational for direct repeats.

Our proof relies on the Cyclic Surgery Theorem, and Theorem A.1
below. In particular, it utilizes the products $b(0,1)$ and
$b(2,1)$. In contrast, the work of \cite{ES1} requires the
products $b(p,q_1), b(p+2,q_2)$ and $b(p+4,q_3)$ (in this setting,
$b(0,1), b(2,1)$ and $b(4,1)$). Since different protein systems
can yield different products, these two theorems may be useful in
complementary situations.

\begin{thm}
\label{noS2xS1andP3} Let $M$ be an orientable Seifert Fiber Space
with incompressible torus boundary. Then there are no slopes
$\alpha$ and $\beta$ on $\partial M$ such that $M(\alpha)=L(0,1)$
and $M(\beta)=L(2,1)$.
\end{thm}

\nibf{Proof.} We first claim that $M$ must be a Seifert Fiber
Space (abbreviated SFS) over $D^2$ with exactly $2$ exceptional
fibers. Note first that $M$ must have base space $D^2$. If the
base space were orientable and not $D^2$, then fillings yield a
SFS over a genus $\geq 1$ surface, not lens spaces. If the base
space were non-orientable then the fibers over an
orientation-reversing curve on the base space would form a Klein
Bottle (as $M$ is orientable), a contradiction. Secondly, if $M$
had $1$ exceptional fiber $M$ would be a solid torus, hence have
compressible boundary. If $M$ had 3 or more exceptional fibers,
then $M(r)$ would either be the connect sum of 3 or more lens
spaces (if $r$ was a fiber slope) or otherwise a SFS over $S^2$
with at least 3 exceptional fibers. In neither case would $M(r)$
be a lens space. (See Orlik \cite{Orl} for more information on
SFSs.) So $M$ must be a SFS over $D^2$ with 2 exceptional fibers.

Fillings along the fiber slope give the connect sum of two lens
spaces, which can never be a lens space. All other fillings of $M$
induce a Seifert fibration over $S^2$ with at least the same two
exceptional fibers. If the filling is to yield a lens space, then
there can be only two exceptional fibers. Hence the induced
fibration on the filling solid torus must be $1$-longitudinal, as
any other filling gives a SFS over $S^2$ with 3 exceptional
fibers.

So we must look at the possible Seifert fiberings of lens spaces
$L(2k,1)$ for $k\in \{0,1\}$. Assume that $M$ has two fillings,
one yielding $L(0,1)$ and the other $L(2,1)$, and let us study
their induced fibrations. We first examine in the general case the
topological constraints that lens spaces place on their possible
fibrations, and then return to the specific fillings above. Note
from above that the induced fibrations would all have the same
exceptional fiber multiplicities as $M$.

Suppose the two exceptional fibers of $M$ have neighborhoods $N_i$
fibered by fibers of slope $\alpha_i={p_i}/{q_i}$ respectively,
where $(p_i,q_i)=1$ and $p_i >1$ and $q_i\neq 0$ for $i=1$, $2$.
The multiplicity of each exceptional fiber is $p_i$. Let
$\beta_1={r_1}/{s_1}$ be a slope on $\partial N_1$ that intersects
$\alpha_1$ exactly once (so that together they form a basis for
$H_1(\partial N_1)$). In this case $|p_1 s_1-q_1 r_1|=1$. Let
$\wt{R}(t_k)=L(2k, 1)$ have the Seifert fibering induced by the
filling slope $t_k$. Then there are Heegaard solid tori $W_i$
which are fibered by the same slopes $\alpha_i$.

Let $\beta^k_2={r^k_2}/{s^k_2}$ be the image of $\beta_1$ under
the gluing map $h:\partial W_1\rightarrow
\partial W_2$ (so we also have $\alpha_2=h(\alpha_1)$). Note that
$\beta^k_2$ depends on the filling $t_k$. We have that $|p_2
s^k_2-q_2 r^k_2|=1$. The meridian $\mu_1=0/1$ of $W_1$ can be
written as a linear combination $m({p_1}/{q_1})+n({r_1}/{s_1})$,
which gives us two equations: $m p_1+n r_1=0$ and $m q_1+n s_1=1$.
By the second equation, $m$ and $n$ are relatively prime, and
since the first equation gives $m p_1= -n r_1$, we have that
$m=-r_1$ and $n=p_1$. Thus
$h(\mu_1)=h(-r_1\alpha_1+p_1\beta_1)=-r_1\alpha_2+p_1 \beta^k_2$.
For the lens space $L(2k,1)$ we then have that:
\begin{eqnarray*}
-r_1 p_2+p_1 r^k_2=2k\\
-r_1 q_2+p_1 s^k_2=1.
\end{eqnarray*}

In the case $L(0,1)$ ($=S^2\times S^1$) we get that $r_1 p_2=p_1
r^0_2$ and as $(p_i,r_i)=1$ this gives us $p_1=p_2$ (so we will
write $p$ for $p_1$ and $p_2$). Thus $r_1=r^0_2$. In other words
the fibration must be the same (with an appropriate choice of
longitudes) on both solid tori.

In the case $L(2,1)$ since the multiplicity must be the same we
have $p_1=p_2=p$ as above. Therefore we have that $-r_1 p + p
r^1_2 =2$ which implies $p|2$, and as $p>1$, this means $p=2$ and
$r_1=r^1_2 -1$. Since we also have $-r_1 q_2+2 s^1_2=1$, we get
$(2 s^1_2-q_2 r^1_2)+q_2= 1$. As $(2 s^1_2-q_2 r^1_2)=\pm 1$ this
implies $q_2=0$ or $\pm 2$, which is impossible as $q_2\neq 0$ and
$(2, q_2)=1$. Hence it is impossible to have a SFS such that one
filling yields $L(0,1)$ and the other $L(2,1)$ \hfill 
\smallskip

\begin{thm}
\label{Rrat}
 $R$ must be rational in the case of direct repeats.
\end{thm}

\nibf{Proof.} Assume $R$ is not rational, so $\wt{R}$ is a
manifold with incompressible torus boundary. Since $O^k$ is
rational for all $k$ by the preceding theorem, $\wt{O^k}$ is a
solid torus, $V_k$, for each $k$. Hence $V_{O^0} \cup \wt{R} =
L(0,1)$, $V_{O^1} \cup \wt{R} = L(2,1)$, $V_{O^2} \cup \wt{R} =
L(4,1)$, and $V_{O^3} \cup \wt{R} = L(6,1)$.

If $\wt{R}$ is not a Seifert Fiber Space, then we have a
contradiction from the Cyclic Surgery Theorem \cite{CGLS1}. In
this case if $\wt{R}(r)$ and $\wt{R}(s)$ are cyclic then $\Delta
(r,s) = 1$, and hence there are at most three such slopes. This
contradicts the fact that we have four fillings which are cyclic.

If $\wt{R}$ is a SFS, then by Theorem~\ref{noS2xS1andP3} it is
impossible to have a SFS such that one filling yields $L(0,1)$ and
the other $L(2,1)$ and so we get a contradiction since $V_{O^0}
\cup \wt{R} = L(0,1)$, $V_{O^1} \cup \wt{R} = L(2,1)$. $\wt{R}$
cannot be a SFS and thus, $\wt{R}$ must be a solid torus, and so
$R$ is rational. \hfill $\Box$

\end{document}